\documentclass[11pt,titlepage]{article}
\usepackage[margin=1.2in]{geometry}
\usepackage{amsmath, amsbsy, amsthm, array, mathrsfs}
\usepackage{graphics}
\usepackage{physics}
\usepackage{epsfig, amssymb,latexsym,verbatim}
\usepackage{graphicx}
\usepackage{color}
\usepackage{dsfont, bbm}
\usepackage{relsize}
\usepackage{subcaption}

\newcommand{\R}{\mathbb{R}}

\newcommand{\noin}{\noindent}
\newcommand{\bee}{\begin{eqnarray*}}
\newcommand{\ene}{\end{eqnarray*}}
\newcommand{\bec}{\begin{center}}
\newcommand{\enc}{\end{center}}
\newcommand{\be}{\begin{equation}}
\newcommand{\ee}{\end{equation}}

\newcommand{\mc}{\mathcal}

\newcommand{\ep}{\varepsilon}
\newcommand{\mb}{\mathbf}
\newcommand{\bs}{\boldsymbol}
\newcommand{\tb}{\textbf}
\newcommand{\pend}{$\blacksquare$}
\newcommand{\vs}{\vskip 3mm}
\newcommand{\bi}{\begin{itemize}}
\newcommand{\ei}{\end{itemize}}

\begin{document}
\baselineskip 3.2ex

\title{\LARGE Asymptotic normality of the least sum of squares of trimmed residuals estimator} 
\vs
\vs
\author{ {\sc 
Yijun Zuo}\\[2ex]
         {\small {\em  
         Department of Statistics and Probability} }\\[.5ex]
         {\small Michigan State University, East Lansing, MI 48824, USA} \\[2ex]
         {\small 
         zuo@msu.edu}\\[6ex]
     }
 \date{\today}
\maketitle

\vskip 3mm
{\small

\begin{abstract}
To enhance the robustness of the  classic least sum of squares (LS) of residuals estimator, Zuo (2022) introduced
the least sum of squares of trimmed (LST) residuals  estimator. The LST enjoys many desired properties and serves well as a robust alternative to the LS. Its asymptotic properties, including strong  and root-n  consistency, have been established whereas the asymptotic normality is left unaddressed.
This article solves this remained problem.

\bigskip
\par

\noindent{\bf Key words and phrase:}  trimmed squares of residuals, robust regression, asymptotics. 
\bigskip
\par
\noindent {\bf Running title:} Asymptotic normality of the LST 
\end{abstract}
}
\setcounter{page}{1}
\section{Introduction}
The general form
of the classical multiple linear regression model is as follows:
\be 
y = \beta_{00} + \beta_{01} x_1 + \cdots + \beta_{0(p-1)}x_{p-1}  + e=(1,\bs{x}')\bs{\beta}_0+ e, \label{general.model}
\ee 
where ' stands for the transpose of a matrix/vector, $y$ is the dependent variable, $\bs{x}=(x_1, x_2, \cdots, x_{p-1})'$ is the independent variable, and $\bs{\beta}_0=(\beta_{00}, \beta_{01}, \cdots, \beta_{0(p-1)})'$ is the regression coefficients (true  unknown regression parameter). 
 Let $\bs{w}=(1, \bs{x}')'$.
Then $y=\bs{w}'\bs{\beta}_0+e$. 
\vs
Assume that one is  given a sample $\mb{Z}^{(n)}
:=\{(\bs{x}'_i, y_i)', i=1,\cdots, n\}$ from the model, where $\bs{x}_i=(x_{i1}, \cdots, x_{i(p-1)})'$  and  wants to estimate the $\bs{\beta}_0$ .  For a given candidate coefficient vector $\bs{\beta}$, call the difference between $y_i$ (observed) and $(1, \bs{x}'_i)'\bs{\beta}$ (predicted by the model), the ith residual, $r_i$ ($\bs{\beta}$ is suppressed).
That is,  with $\bs{w}_i:=(1,\bs{x}'_i)'$
\be {r}_i=y_i-(1, \bs{x'_i})\bs{\beta}=y_i-\bs{w}'_i\bs{\beta}.\label{residual-1.eqn}
\ee
To estimate $\bs{\beta}_0$, the classic \emph{least squares} (LS) estimator is the minimizer of the sum of squares of residuals
\be\widehat{\bs{\beta}}_{ls}=\arg\min_{\bs{\beta}\in\R^p} \sum_{i=1}^n r^2_i. \label{ls.eqn}
\ee
A straightforward calculus derivation leads to
$$
\widehat{\bs{\beta}}_{ls}=(\bs{X}_n \bs{X}'_n)^{-1}\bs{X}'_n\bs{Y}_n,
$$
where $\bs{Y}_n=(y_1,\cdots, y_n)'$, $\bs{X}_n=(\bs{w}_1,\cdots, \bs{w}_n)'$ and $\bs{x}_1,\cdots, \bs{x}_n$ are assumed to be linearly independent (i.e. $\bs{X}_n$ has a full rank).
\vs
Due to its great computability and optimal properties when the error $e$ follows a normal $\mc{N}(\mu,\sigma^2)$ distribution, the least squares estimator is the  most popular  in practice across multiple disciplines and the benchmark in the multiple linear regression.\vs
It, however, can behave badly when the error distribution is slightly departed from the normal distribution,
particularly when the errors are heavy-tailed or contain outliers.
In fact, both $L_1$ (squared residuals replaced by absolute residuals in (\ref{ls.eqn})) and $L_2$ (LS) estimators have a pathetic $0\%$ asymptotic breakdown point (see Section 3.1 of Zuo (2022)), in sharp contrast to the $50\%$ of the least trimmed residuals (Rousseeuw (1984)).  The latter is one of the most robust alternatives to the least squares estimator.
\vs
Seminal papers by Box (1953)  and Tukey (1960)  were the impetus for robust statistical procedures. The theory of robust statistics blossomed in the 1960s -- 1980s. 
 Robust alternatives to the least squares regression estimator are abundant in the literature.
 The most popular
ones are, among others, M-estimators (Huber(1964)), least median squares (LMS) and least
trimmed squares (LTS) estimators (Rousseeuw (1984)), S-estimators (Rousseeuw and Yohai
(1984)), MM-estimators (Yohai (1987) ), $\tau$ -estimators (Yohai and Zamar (1988) ) and maximum
depth estimators (Rousseeuw and Hubert (1999) and Zuo (2021a,b)).
 \vs
 Among all robust alternatives,  in practice, the LTS is one of the most prevailing crossing multiple disciplines. Its idea is simple, ordering the squared residuals and then trimming the larger ones and keeping at least $\lfloor n/2\rfloor$ squared residuals, where $\lfloor ~\rfloor$ is the floor function, the  minimizer of the sum of those \emph{trimmed squared residuals} is called an LTS estimator:
\be
\widehat{\bs{\beta}}_{lts}:=\arg\min_{\bs{\beta}\in \R^p} \sum_{i=1}^h (r^2)_{(i)}, \label{lts.eqn}
\ee
where $(r^2)_{(1)}\leq (r^2)_{(2)}\leq \cdots \leq (r^2)_{(n)}$ are the ordered squared residuals and $\lfloor n/2\rfloor \leq h \leq n$.
\vs
Realized the high variability of $\widehat{\bs{\beta}}_{lts}$, Zuo (2022) introduced the least sum of squares of trimmed (LST) residuals  estimator. Instead of trimming after squaring of residuals as LTS does, the LST, employing a depth/outlyingness based scheme, trims residuals first then squares the rest. The minimizer of the sum of squares of trimmed residuals is called an LST estimator. Before formally introducing
LST in Section 2, let us first appreciate the difference among the LS, the LTS, and the LST procedures.\vs
\bec
\begin{figure}[h]
    \centering
    \vspace*{-8mm}
    \includegraphics[width=8cm, height=7cm] 
    {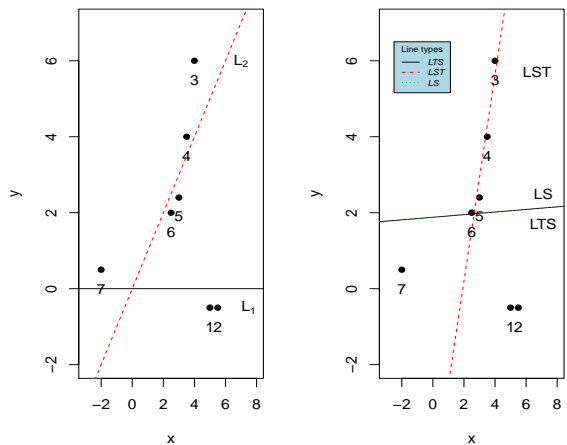}
        \caption{\small Left panel: plot of seven artificial points and two candidate lines ($L_1$ and $L_2$), which line would you  pick?
        Sheerly based on the trimming scheme and objective function value, if one uses the number $h=\lfloor n/2\rfloor+\lfloor (p+1)/2\rfloor$, given on page 132 of  Rousseeuw and Leroy (1987) (RL87), for achieving the best robustness, that is, employing four squared residuals, then LTS prefers $L_1$ to $L_2$ whereas LST reverses the preference.  Right panel:  the same seven points are fitted by LTS,  LST, and  LS (benchmark). The solid black  line is the LTS given by ltsReg. Red dashed line is given by the LST, and green dotted line is given by the LS - which is identical to the LTS line in this case.}
        \label{Fig:one}
\end{figure}
\vspace*{-10mm}
\enc
\vs
\noin
\tb{Example 1.1}
For illustration purpose, we borrow a small data set from Zuo (2022)   and with $x=(5, 5.5, 4, 3.5, 3, 2.5, -2)$ and $y=(-.5, -.5, 6, 4, 2.4, 2, 0.5)$. That is, sample size $n=7$ and dimension $p=2$. The data are plotted in the left panel  of Figure \ref{Fig:one}. We also provide two candidate regression lines $L_1$ ($y=0$) and $L_2$ ($y=x$). Which one would you pick to represent the overall pattern of the data set?
 \vs
Intuitively points $1$ and $2$ are outliers and $L_2$ should be preferred in the left panel  of Figure \ref{Fig:one}. But that is not the case if one employs the R function ltsReg for the LTS which gives the line in the right panel of Figure \ref{Fig:one}, along with it is the benchmark LS line (identical to the LTS-ltsReg line) and the line by  the procedure LST (see Zuo (2022) ). Obviously, both the LS and the LTS-ltsReg lines, influenced by the two outliers,  miss the overall linear pattern of the data whereas the LST line resists the outlier's influence and still catches the overall linear pattern.
\hfill \pend
\vs
Similar examples for an increased sample size or real data set are given in Zuo (2022).\vs
The rest of article is organized as follows. Section 2 formally introduces the least sum of squares of trimmed residuals estimator and establishes or summarizes its preliminary properties which will be useful in Section 3 where its asymptotic normality is established via stochastic euqicontinuity. Inference procedures based on the asymptotic normality and bootstrapping are addressed in Section 4. Concluding remarks end the article in Section 5.

\section{The LST and its preliminary properties}
\subsection{Depth trimming and LST} The LTS squares residuals first and then trims,
the LST, on the other hand, employing a depth or outlyingness based scheme,  trims residuals first then squares the left. 
\vs 
The LTS utilizes a rank-based trimming scheme. The latter usually focuses only on the relative position of points with respect to others
and ignores the magnitude of the point and the relative distance between points. Zuo (2006) 
and Wu and Zuo (2009) discussed an alternative trimming scheme, which exactly catches
these two important attributes (magnitude and relative distance). It orders data from a
center (the median) outward and trims the points that are far away from the center. This is
known as depth-based trimming. In other words, the depth-based trimming
scheme trims points that lie on the outskirts (i.e. points that are less deep, or outlying). The depth (or,
equivalently, outlyingness ) of a point x is defined to be
\be
D\big(x, x^{(n)}\big)=\frac{|x-\mbox{Med}(x^{(n)})|}{\mbox{MAD}(x^{(n)})},  \label{outlyingness.eqn}
\ee
where $x^{(n)}=\{x_1, \cdots, x_n\}$ is a data set in $\R^1$,  Med$(x^{(n)})=\mbox{median}(x^{(n)})$ is the median of the data points, and  MAD$(x^{(n)})=\mbox{Med}(\{|x_i-\mbox{Med}(x^{(n)})|,~ i=1,2, \cdots, n\})$ is the median of absolute deviations to the center (median). 
For a given data set $\mb{Z}^{(n)}=\{(\bs{x}'_i, y_i)'\}$ in $\R^{p}$ and a $\bs{\beta} \in \R^p$, define
\begin{align}
\mu_n(\bs{\beta}):=\mu(\bs{Z}^{(n)},\bs{\beta})&=\mbox{Med}_i\{r_i\}, \label{med.eqn}\\
\sigma_n(\bs{\beta}):=\sigma(\bs{Z}^{(n)},\bs{\beta})&=\mbox{MAD}_i\{r_i\}, \label{mad.eqn}
\end{align}
where operators Med and MAD are used for discrete data sets (and distributions as well) and $r_i$ defined in (\ref{residual-1.eqn}).   
For a given $\alpha$ (throughout constant $\alpha\geq 1$, default value is one) in the depth trimming scheme, consider the quantity
\be
Q^n(\bs{\beta}):=Q(\bs{Z}^{(n)}, \bs{\beta}, \alpha)= \frac{1}{n}\sum_{i=1}^{n}r_i^2\mathds{1}\left( \frac{|r_i-\mu(\bs{Z}^{(n)},\bs{\beta})|}{\sigma(\bs{Z}^{(n)},\bs{\beta})}\leq \alpha\right),\label{objective.eqn}
\ee
 where $\mathds{1}(A)$ is the indicator of $A$ (i.e., it is one if A holds and zero otherwise).
Namely, residuals with their outlyingness (or depth) greater than $\alpha$ will be trimmed.
When there is a majority ($\geq \lfloor(n+1)/2\rfloor$) identical $r_i$s, we define $\sigma(\mb{Z}^{(n)}, \bs{\beta})=1$.
Minimizing $Q(\bs{Z}^{(n)}, \bs{\beta}, \alpha)$, one gets the \emph{least} sum of \emph{squares} of {\it trimmed} (LST)  residuals  estimator,
\be
\widehat{\bs{\beta}}^n_{lst}:=\widehat{\bs{\beta}}_{lst}(\mb{Z}^{(n)}, \alpha)=\arg\min_{\bs{\beta}\in \R^p}Q(\bs{Z}^{(n)}, \bs{\beta}, \alpha).\label{lst.eqn}
\ee
One might take it for granted that the  minimizer of  $Q(\bs{Z}^{(n)}, \bs{\beta}, \alpha)$ always exists.
Does the right-hand side (RHS) of (\ref{lst.eqn}) always have a minimizer? If yes, will it be unique?\vs
\noin
\subsection{Existence, Uniqueness, Fisher consistency, and Equivariance}
\vs
\noindent
\tb{Existence and Uniqueness}~~
 For simplicity of description, we write $\mathds{1}(i, \bs{Z}^{(n)}, \bs{\beta}, \alpha)$ for $\mathds{1}\left( {|r_i-\mu_n(\bs{\beta})|}/{\sigma_n(\bs{\beta})}\leq \alpha\right)$.  \vs
\vs
\noindent
\tb{Theorem 2.1} We have
\bi
\item[(i)] $\widehat{\bs{\beta}}^n_{lst}$ always exist;
\vs
\item[(ii)] $\widehat{\bs{\beta}}^n_{lst}$ is unique if
$\bs{M}:=\bs{M}(\bs{Y}_n, \bs{X}_n, \bs{\beta}, \alpha)=\sum_{i=1}^n \bs{w}_i\bs{w}'_i\mathds{1} (i, \bs{Z}^{(n)}, \bs{\beta}, \alpha)$ is invertible.

\ei
\vs
\noindent
\tb{Proof}: (ii) was covered by Theorem 2.3 of Zuo (2022) . (i) was also proved in Theorem 2.2 there with \emph{an extra assumption} though.
We now show (i) without any assumption. 
By the proof of Theorem 2.1 of Zuo (2022) , it is seen that
\[ 
\frac{\partial Q^n(\bs{\beta})}{\partial \bs{\beta}} =\frac{2}{n}\sum_{i=1}^n  r_i \bs{w}_i \mathds{1} (i, \bs{Z}^{(n)}, \bs{\beta}, \alpha).
\]
Furthermore,
\[\frac{\partial^2 Q^n(\bs{\beta})}{\partial \bs{\beta}^2}=\frac{2}{n}\sum_{i=1}^n \bs{w}_i\bs{w}'_i\mathds{1} (i, \bs{Z}^{(n)}, \bs{\beta}, \alpha).
\]
Call the matrix on the RHS above as H (Hessian matrix). It is readily seen that $H$ is positive semidefinite. Hence, $Q^n(\bs{\beta})$ is convex  and twice continuously differentiable in $\bs{\beta}$. Consequently the global minimum of $Q^n(\bs{\beta})$, $\widehat{\bs{\beta}}^n_{lst}$, always exists \hfill \pend
\vs
\noindent
\tb{Remark 2.1}
\bi
\item[] Uniqueness is indispensable for later asymptotic normality establishment.
A sufficient condition for $M$ in (ii) being invertible is that $\bs{x}_1, \cdots, \bs{x}_n$ are linearly independent, or the $\bs{X}_n$ has a full rank. \hfill \pend
\ei
\vs
There is a counterpart of Theorem 2.1 at the population setting. To that end, we first need to have the counterparts of (\ref{objective.eqn}) and (\ref{lst.eqn}) at the population setting. \vs
Throughout $F_{\mb{z}}$ always stands for the distribution of random vector $\mb{z}$ unless otherwise stated.
Write $F_{(\bs{x'}, y)}$ for the joint distribution of $\bs{x}'$ and $y$ in the model (\ref{general.model}) and  $\bs{w}=(1,\bs{x}')'$
\begin{align}
r:=&r(F_{(\bs{x'}, y)}, \bs{\beta})=y-\bs{w}'\bs{\beta}, \label{residual.eqn}\\
  \mu(F_r):=&\mu(F_{(\bs{x'}, y)}, \bs{\beta})=\mbox{Med}(F_r),\label{popultaion.med.eqn}\\
  \sigma(F_r):=&\sigma(F_{(\bs{x'}, y)}, \bs{\beta})=\mbox{MAD}(F_r),\label{population.mad.eqn}
  \end{align}
hereafter we assume that $\mu$ and $\sigma$ exist uniquely.
The population counterparts of (\ref{objective.eqn}) and (\ref{lst.eqn}) are respectively:
\begin{align}
Q(F_{(\bs{x'}, y)},\bs{\beta}):=Q(F_{(\bs{x'}, y)},\bs{\beta}, \alpha):&=\int(y-\bs{w}'\bs{\beta})^2\mathds{1}\left( \frac{|y-\bs{w}'\bs{\beta}-\mu(F_r)|}{\sigma(F_r)}\leq \alpha\right)dF_{(\bs{x'}, y)}, \label{Q.eqn}\\[1ex]
\bs{\beta}_{lst}:=\bs{\beta}_{lst}(F_{(\bs{x'}, y)}):&=\bs{\beta}_{lst}(F_{(\bs{x'}, y)}, \alpha) :=\arg\min_{\bs{\beta}\in \R^p}Q(F_{(\bs{x'}, y)}, \bs{\beta}, \alpha). \label{lst-def.eqn}
\end{align}
\vs
\noindent
\tb{Theorem 2.2}
\bi
\item[(i)] $\bs{\beta}_{lst}(F_{(\bs{x'}, y)})$ always exists;
\vs
\item[(ii)]  $\bs{\beta}_{lst}(F_{(\bs{x'}, y)})$ is unique if $E_{(\bs{x}', y)}\left(\bs{w}\bs{w}'\mathds{1} (F_{(\bs{x'}, y)}, \bs{\beta}, \alpha)\right)$ is invertible, where $\mathds{1} (F_{(\bs{x'}, y)}, \bs{\beta}, \alpha)$ stands for $\mathds{1}\left( {|y-\bs{w}'\bs{\beta}-\mu(F_r)|}/{\sigma(F_r)}\leq \alpha\right)$.
\ei
\vs
\noindent
\tb{Proof:} This is analogue to that of Theorem 2.1.
\vs
(i) Take the first-order and second-order derivative of $Q(F_{(\bs{x}', y)},\bs{\beta})$ with respect to $\bs{\beta}$, in light of Lebesgue dominating theorem, we have
\begin{align}
\frac{\partial Q(F_{(\bs{x}', y)},\bs{\beta})}{\partial \bs{\beta}}&=2 E_{(\bs{x}', y)}\big(r\bs{w}\mathds{1}(F_{(\bs{x'}, y)}, \bs{\beta}, \alpha)\big),\\[1ex]
\frac{\partial^2 Q(F_{(\bs{x}', y)},\bs{\beta})}{\partial \bs{\beta}^2}&=2 E_{(\bs{x}', y)}\big(\bs{w}\bs{w}'\mathds{1}(F_{(\bs{x'}, y)}, \bs{\beta}, \alpha) \big).
\end{align}
It is readily seen that the RHS matrix in the last equation is positive semidefinite, hence $Q(F_{(\bs{x}', y)},\bs{\beta})$ is twice continuously differentiable and convex in $\bs{\beta}$. Consequently, the global minimum of $Q(F_{(\bs{x}', y)},\bs{\beta})$, $\bs{\beta}_{lst}(F_{(\bs{x'}, y)})$, always exists.
\vs
(ii) When $E_{(\bs{x}', y)}\left(\bs{w}\bs{w}'\mathds{1} (F_{(\bs{x'}, y)}, \bs{\beta}, \alpha)\right)$ is invertible then $Q(F_{(\bs{x}', y)},\bs{\beta})$ is strictly convex in $\bs{\beta}$, the uniqueness follows.
\hfill \pend
\vs
\noindent
\tb{Remark 2.2}
\bi
\item[] Existence and uniqueness  are also established in Zuo (2022)  with much more assumptions whereas here we do it without any assumption or with one assumption, respectively.
\hfill \pend
\ei
\vs
 \noindent
 \tb{Fisher consistency}~~ Next
 we like to show that  $\bs{\beta}_{lst}(F_{(\bs{x'}, y)})$ is identical to the true unknown parameter $\bs{\beta}_0$ under some assumptions - which is called Fisher consistency of the estimation functional. Recall our general model:
$
y=\bs{w}'\bs{\beta}_0+e. 
$
\vs
\noin
 \tb{Theorem 2.3}
  \bi
  \item[] \hspace*{-5mm}$\bs{\beta}_{lts}(F_{(\bs{x}', y)})=\bs{\beta}_0$  provided that

 \item[(i)] $E_{(\bs{x}', y)}\left(\bs{w}\bs{w}' \mathds{1} (F_{(\bs{x'}, y)}, \bs{\beta}, \alpha)\right)$ is invertible,

 \item[(ii)] $E_{(\bs{x}', y)}
 \left( e\bs{w} \mathds{1}\left( \frac{|e-\mu({F_e})|}{\sigma(F_e)}\leq \alpha\right)\right)=\mb{0}$.
 \ei
 \vs
 \noindent
 \tb{Proof} By theorem 2.2, (i) guarantees the unique existence  of $\bs{\beta}_{lts}(F_{(\bs{x}', y)})$ which is the unique solution of the system of the equations
\[
 \int(y-\bs{w}'\bs{\beta})\bs{w}\mathds{1} (F_{(\bs{x'}, y)}, \bs{\beta}, \alpha) dF_{(\bs{x}', y)}(\bs{x}, y)=\mb{0}.
 \]
 Notice that $y-\bs{w}'\bs{\beta}=-\bs{w}'(\bs{\beta}-\bs{\beta}_0)+e$, insert this into the above equation we have
 \[
 \int(-\bs{w}'(\bs{\beta}-\bs{\beta}_0)+e)\bs{w}\mathds{1}\left(\frac{|-\bs{w}'(\bs{\beta}-\bs{\beta}_0)+e-\mu(F_r)|}{\sigma(F_r)}\right) dF_{(\bs{x}', y)}(\bs{x}, y)=\mb{0}.
 \]
 By (ii) it is readily seen that $\bs{\beta}=\bs{\beta}_0$ is a solution of the above system of equations. Uniqueness leads to the desired result.
 \hfill \pend
\vs
\noindent
\tb{Remark 2.3}
\bi
\item[] Fisher consistency is also proved in Zuo (2022)  under four assumptions though. \hfill \pend
\ei
\vs
\noin
\tb{Equivariance}
A regression functional $\bs{T}(\cdot)$ is \emph{regression}, \emph{scale}, and \emph{affine} equivariant, (see Zuo (2021a))  if, respectively,
\begin{align*}
\bs{\beta}^*(F_{(\bs{w},~ y+~\bs{w}'\bs{b})})=\bs{\beta}^*(F_{(\bs{w},~ y)})+\bs{b}, \forall ~\bs{b}\in \R^p;\\
\bs{\beta}^*(F_{(\bs{w}, ~sy)})=s \bs{\beta}^*(F_{(\bs{w}, ~y)}), \forall ~s \in \R;\\
\bs{\beta}^*(F_{(\bs{A}'\bs{w},~ y)})=\bs{A}^{-1}\bs{\beta}^*(F_{(\bs{w},~ y)}),~ \forall \mbox{~ nonsingular ${p\times p}$ matrix  $\bs{A}$}. 
\end{align*}
Namely, $\bs{T}(\cdot)$ does not depend on the underlying coordinate system and measurement scale.\vs
For definition of \emph{regression, scale, and affine equivariance} of a regression estimator  at sample setting, see Zuo (2022) .
\vs
\noindent
\tb{Theorem 2.4}~ $\widehat{\bs{\beta}}^n_{lst}$ and $\bs{\beta}_{lst}$ are regression, scale, and affine equivariant at sample and population settings, respectively. \hfill \pend
\vs
\noindent
\section{Asymptotic normality of the LST}

For a given sample $\bs{Z}^{(n)}=\{\bs{Z}_i\}=\{(\bs{x}'_i, y_i)\}$, $i=1,2,\cdots, n$, write  $F^n_\mb{Z}$ as the sample version of $F_\mb{Z}:=F_{(\bs{x'}, y)}$ based on $\mb{Z}^{(n)}$. It will be used interchangeably with $P_n$ or $\bs{Z}^{(n)}$. 

\noindent
\subsection{Strong consistency}
To show that $\widehat{\bs{\beta}}^n_{lst}$ converges to $\bs{\beta}_{lst}$ almost surely, one can take the approach given in Section 4.2 of Zuo (2022) . But here we take a different directly approach.\vs
Following the notations of Pollard (1984) (P84) , write
\begin{align*}
Q(\bs{\beta}, P):&=Q(F_{\bs{Z}}, \bs{\beta}, \alpha)=P [(y-\bs{w}'\bs{\beta})^2\mathds{1}(F_{(\bs{x}',y)},\bs{\beta}, \alpha)]=Pf,\\[1ex]
 Q(\bs{\beta}, P_n):&=Q(F^n_{\bs{Z}}, \bs{\beta}, \alpha)=\frac{1}{n}\sum_{i=1}^n r^2_i\mathds{1}(i, \bs{Z}^{(n)}, \bs{\beta}, \alpha)= P_nf,
\end{align*}
where $f:=f(\bs{x}, y, \bs{\beta}, \alpha)=(y-\bs{w}'\bs{\beta})^2\mathds{1}(F_{(\bs{x}',y)},\bs{\beta}, \alpha)$. 
\vs
Under corresponding assumptions in Theorems 2.1 and 2.2, $\widehat{\bs{\beta}}^n_{lst}$ and $\bs{\beta}_{lst}$ are unique minimizers of $Q(\bs{\beta}, P_n)$ and $Q(\bs{\beta}, P)$ over $\bs{\beta}\in \R^p$, respectively.\vs

To show that $\widehat{\bs{\beta}}^n_{lst}$ converges to $\bs{\beta}_{lst}$ almost surely, it suffices to prove that $Q(\widehat{\bs{\beta}}^n_{lst}, P) \to Q(\bs{\beta}_{lst}, P)$ almost surely, because $Q(\bs{\beta}, P)$ is bounded away from $Q(\bs{\beta}_{lst}, P)$ outside each neighborhood of $\bs{\beta}_{lst}$ in light of continuity and compactness (also see Lemma 4.3 of Zuo (2022)).\vs

By theorems 2.1 and 2.2, assume, without loss of generality (w.l.o.g.), that  $\widehat{\bs{\beta}}^n_{lst}$ and $\bs{\beta}_{lst}$ belong to a ball centered at $\bs{\beta}_{lst}$ with large enough radius $r_0$, $B(\bs{\beta}_{lst}, r_0)$ (see Section 4.2 of Zuo (2022)).
 Assume, w.l.o.g., that $\Theta=B(\bs{\beta}_{lst}, r_0)$ is our parameter space of $\bs{\beta}$ hereafter. Define a class of functions for a fixed $\alpha\geq 1$
\[
\mathscr{F} (\bs{\beta})=\left\{f(\bs{x}, y, \bs{\beta}, \alpha)= (y-\bs{w}'\bs{\beta})^2\mathds{1}(F_{(\bs{x}',y)},\bs{\beta}, \alpha): \bs{\beta} \in \Theta\right\}.
\]
\vs
If we prove uniform almost sure convergence of $P_n$ to $P$ over $\mathscr{F}$ (generalized Glivenko-Cantelli theorem, see Lemma 3.1 below), then we can deduce almost surely that $Q(\widehat{\bs{\beta}}^n_{lst}, P) \to Q(\bs{\beta}_{lst}, P)$ from
\begin{align*}
Q(\widehat{\bs{\beta}}^n_{lst}, P_n)-Q(\widehat{\bs{\beta}}^n_{lst}, P) &\to 0 ~~(\mbox{in light of Lemma 3.1), ~~and}\\[1ex]
Q((\widehat{\bs{\beta}}^n_{lst}, P_n)\leq Q(\bs{\beta}_{lst}, P_n)& \to Q(\bs{\beta}_{lst}, P)\leq Q(\widehat{\bs{\beta}}^n_{lst}, P). 
\end{align*}

Above discussions and arguments have led to
\vs
\noindent
\tb{Theorem 3.1}. Under corresponding assumptions in Theorems 2.1 and 2.2 for uniqueness of $\widehat{\bs{\beta}}^n_{lst}$ and $\bs{\beta}_{lst}$ respectively, we have
$\widehat{\bs{\beta}}^n_{lst}$ converges almost surely to $\bs{\beta}_{lst}$ (i.e. $\widehat{\bs{\beta}}^n_{lst}-\bs{\beta}_{lst}=o(1)$, a.s.). \hfill \pend
\vs
\noindent
\tb{Lemma 3.1} [Zuo(2022)]. $\sup_{f\in \mathscr{F}}|P_n f-P f| \to 0$ almost surely.  \hfill \pend
\vs
\noin
\subsection{Asymptotic normality}

Instead of treating the root-n consistency separately as Zuo (2022)  did, we will establish asymptotic normality of $\widehat{\bs{\beta}}^n_{lst}$ directly via stochastic equicontinuity (see page 139 of P84, or the supplementary of Zuo (2020)), and consequently obtain  the root-n consistency of $\widehat{\bs{\beta}}^n_{lst}$ 
as a by-product of the asymptotic normality.

\vs
\emph{Stochastic equicontinuity} refers to a sequence of
stochastic processes $\{Z_n(t): t \in T\}$ whose shared index set $T$ comes equipped
with a semi metric $d(\cdot, \cdot)$. (a semi metric has all
the properties of a metric except that $d(s, t) = 0$ need not imply that $s$ equals
$t$.)
\vs
\noin
\tb{Definition 3.1} [IIV. 1, Def. 2 of P84].  Call ${Z_n}$ stochastically equicontinuous at $t_0$  if for each $\eta > 0$
and $\epsilon > 0$ there exists a neighborhood $U$ of $t_0$ for which
\be
\limsup P\left(\sup_{U} |Z_n(t) - Z_n(t_0) | > \eta\right) < \epsilon.  \label{se.eqn}
\ee
\hfill~~~~~~~~~~~~~\pend
\vs
Because stochastic equicontinuity bounds $Z_n$ uniformly over the neighborhood
$U$, it also applies to any randomly chosen point in the neighborhood. 
If ${\tau_n}$ is a sequence of random elements of $T$ that converges in probability
to $t_0$, then
\be
Z_n(\tau_n)-Z_n(t_0)\to 0\mbox{~in probability,}
\ee
because, with probability tending to one, $\tau_n$ will belong to each $U$.
The form above will be easier to apply, especially when behavior of a particular ${\tau_n}$ sequence is under investigation.
\vs
Again following the notations of P84. Suppose $\mathscr{F} = \{ f(\cdot, t): t\in T\}$, with $T$ a subset of $\R^k$, is a collection of
real, P-integrable functions on the set $S$ where $P$ (probability measure) lives.
Denote by $P_n$ the
empirical measure formed from $n$ independent observations on $P$, and define the empirical process $E_n$ as the signed measure $n^{1/2}(P_n - P)$. Define
\begin{align*}
F(t) &= P f(\cdot, t),\\
F_n(t) &= P_n f(\cdot, t).
\end{align*}
Suppose $f(\cdot, t)$ has a linear approximation near the $t_0$ at which $F(\cdot)$ takes
on its minimum value:
\be
 f(\cdot, t) = f(\cdot, t_0) + (t - t_0)'\nabla(\cdot) + |t - t_0|r(\cdot, t). \label{taylor.eqn}
\ee
For completeness set $r(\cdot, t_0) = 0$, where $\nabla$ (differential operator) is a vector of $k$ real functions on
$S$. We cite theorem 5 of IIV.1 of P84  (page 141) for the asymptotic normality of $\tau_n$.
\vs
\noindent
\tb{Lemma 3.2} .  Suppose $\{\tau_n\}$ is a sequence of random vectors converging in
probability to the value $t_0$ at which $F(\cdot)$ has its minimum. Define $r(\cdot, t)$ and the
vector of functions $\nabla(\cdot)$ by (\ref{taylor.eqn}). If
\bi
\item[(i)] $t_0$ is an interior point of the parameter set $T$; \vspace*{-2mm}
\item[(ii)] $F(\cdot)$ has a non-singular second derivative matrix $V$ at $t_0$;\vspace*{-2mm}
\item[(iii)] $F_n(\tau_n) = o_p(n^{-1}) + \inf_{t}F_n(t)$;\vspace*{-2mm}
\item[(iv)] the components of $\nabla(\cdot)$ all belong to $\mathscr{L}^2(P)$;\vspace*{-2mm}
\item[(v)] the sequence $\{E_{n}(\cdot, t)\}$ is stochastically equicontinuous at $t_0$ ;\vspace*{-2mm}
\ei
then
\[n^{1/2}(\tau_n - t_0) \stackrel{d} \longrightarrow  {\cal{N}}(O, V^{-1}[P(\nabla\nabla') - (P\nabla)(P\nabla)']V^{-1}).
\]
\vs
In order to apply the Lemma, we first realize that in our case, $\widehat{\bs{\beta}}^n_{lst}$ and $\bs{\beta}_{lst}$ correspond to $\tau_n$ and $t_0$ (assume, w.l.o.g. that $\bs{\beta}_{lts}=\mb{0}$ in light of regression equivariance); $\bs{\beta}$ and $\Theta$ correspond to $t$ and $T$; $f(\cdot, t):= f(\cdot, \cdot, \bs{\beta}, \alpha)$ and $\alpha$ is a fixed constant. 
In our case,
\[
\nabla(\bs{x}, y, \bs{\beta}, \alpha)=\frac{\partial}{\partial \bs{\beta}} f(\bs{x}, y, \bs{\beta}, \alpha)=2(y-\bs{w}'\bs{\beta})\bs{w}\mathds{1}(F_{(\bs{x}', y)}, \bs{\beta},\alpha).
\]
We will have to assume that $P(\nabla^2_i)=P(4(y-\bs{w}'\bs{\beta})^2{w}^2_i\mathds{1}(F_{(\bs{x}', y)}, \bs{\beta},\alpha)$ exists to meet (iv) of the lemma, where $i\in \{1,\cdots, p\}$ and $\bs{w}'=(w_1, \cdots, w_p)=(1, \bs{x}')$. It is readily seen that a sufficient condition for this assumption to hold is the existence of $P(x^2_i)$. In our case,
$V=2 P(\bs{w}\bs{w}'\mathds{1}(F_{(\bs{x}', y)}, \bs{\beta},\alpha)$, we will have to assume that it is invertible when $\bs{\beta}$ is replaced by $\bs{\beta}_{lst}$ (it is covered by the assumption in Theorem 2.2)  to meet (ii) of the lemma. In our case,
\[r(\cdot, t)=\left(\frac{\bs{\beta}'}{\|\bs{\beta}\|}V/2 \frac{\bs{\beta}}{\|\bs{\beta}\|} \right)\|\bs{\beta}\|.\]
We will assume that $\lambda_{min}$ and $\lambda_{max}$ are the minimum and maximum eigenvalues of positive semidefinite matrix $V$ overall $\bs{\beta}\in \Theta$ and a fixed $\alpha \geq 1$.
\vs
\noindent
\tb{Theorem 3.2} Assume that\vspace*{-2mm}
\bi
\item[(i)] the uniqueness assumptions for $\widehat{\bs{\beta}}^n_{lst}$ and $\bs{\beta}_{lst}$ in theorems 2.1 and 2.2 hold respectively;\vspace*{-2mm}
\item[(ii)] $P({x^2_i})$ exists; \vspace*{-2mm}
\ei
then
\[n^{1/2}(\widehat{\bs{\beta}}^n_{lst} - \bs{\beta}_{lst}) \stackrel{d} \longrightarrow  {\cal{N}}(O, V^{-1}[P(\nabla\nabla') - (P\nabla)(P\nabla)']V^{-1}),
\]
where $\bs{\beta}$ in $V$ and $\nabla$ is replaced by $\bs{\beta}_{lst}$ (which could be assumed to be zero).
\vs
\noindent
\tb{Proof}: To apply Lemma 2.5, we need to verify the five conditions, among them only (iii) and (v) need to be addressed, all others are satisfied trivially. For (iii), it holds automatically since our $\tau_n=\widehat{\bs{\beta}}^n_{lst}$ is defined to be the minimizer of  $F_n(t)$ over $t\in T (=\Theta)$.\vs
So the only condition that needs to be verified is the (v), the stochastic equicontinuity of $\{E_nr(\cdot, t)\}$ at $t_0$. For that, we will appeal to the Equicontinuity Lemma (VII.4 of P84, page 150).
To apply the Lemma, we will verify that the condition for the random covering numbers satisfy the uniformity condition. To that end, we look at the class of functions for a fixed $\alpha\geq 1$
\[
\mathscr{R}(\bs{\beta})=\left\{r(\cdot, \cdot, \alpha, \bs{\beta})=\left(\frac{\bs{\beta}'}{\|\bs{\beta}\|}V/2 \frac{\bs{\beta}}{\|\bs{\beta}\|} \right)\|\bs{\beta}\|:~ \bs{\beta}\in \Theta \right\}.
\]
Obviously, $\lambda_{max} r_0/2$ is an envelope for the class $\mathscr{R}$ in $\mathscr{L}^2(P)$, where $r_0$ is the radius of the ball $\Theta=B(\bs{\beta}_{lts}, r_0)$. We now show that the covering numbers of $\mathscr{R}$ are uniformly bounded, which amply suffices for the Equicontinuity Lemma.  For this, we will invoke Lemmas II.25 and II.36 of P84. 
 To apply Lemma II.25, we need to show that the graphs of functions in $\mathscr{R}$ have only polynomial discrimination.
 \vs
The graph of a real-valued function $f$ on a set $S$ is defined as the subset (see page 27 of P84 )
$$G_f = \{(s, t): 0\leq t \leq f(s) ~\mbox{or}~ f(s)\leq t \leq 0, s \in S \}.$$
\vs
  The graph of $r(\bs{x}, y, \alpha, \bs{\beta})$ contains a point $(\bs{x}, y, t)$, $t\geq 0$
  if and only if $\left(\frac{\bs{\beta}'}{\|\bs{\beta}\|}V/2 \frac{\bs{\beta}}{\|\bs{\beta}\|} \right)\|\bs{\beta}\| \geq t$ for all $\bs{\beta} \in \Theta$. 
 Equivalently, the graph of $r(\bs{x}, y, \alpha, \bs{\beta})$ contains a point $(\bs{x}, y, t)$, $t\geq 0$ if and only if $\lambda_{min}/2\|\bs{\beta}\|\geq t$.
For a collection of $n$ points $(\bs{x}'_i, y_i, t_i)$ with $t_i\geq 0$, the graph picks out those points satisfying $\lambda_{min}/2\|\bs{\beta}\|- t_i\geq 0$. Construct from
$(\bs{x}_i, y_i, t_i)$ a point $z_i=t_i$ in $\R$. On $\R$ define a vector space $\mathscr{G}$ of functions
\[
g_{a, b}(x)=ax+b,~~ a,~ b \in \R.
\]
By Lemma 18 of P84, the sets $\{g\geq 0\}$, for $g \in \mathscr{G}$, pick out only a polynomial number of subsets from $\{z_i\}$; those sets corresponding to functions in $\mathscr{G}$ with $a=-1$ and $b=\lambda_{min}/2\|\bs{\beta}\|$ pick out even fewer subsets from  $\{z_i\}$. Thus the graphs of functions in $\mathscr{R}$ have only    polynomial discrimination.
\hfill \pend
\vs
\noindent
\section{Inference procedures}
In order to utilize the asymptotic normality result in Theorem 3.2, we need to figure out the asymptotic covariance. Assume that $\bs{z}=(\bs{x}', y)'$ follows  elliptical distributions $E(g; \bs{\mu}, \bs{\Sigma})$ with density
$$
f_{\bs{z}}(\bs{x}', y)=\frac{g(((\bs{x}',y)'-\bs{\mu})'\bs{\Sigma}^{-1}((\bs{x}',y)'-\bs{\mu}))}{\sqrt{\det(\bs{\Sigma})}},
$$ 
where 
$\bs{\mu}\in \R^p$ and $\bs{\Sigma}$  a positive definite matrix of size $p$ which is proportional to the covariance matrix if the latter exists. We assume the function $g$ to have a strictly negative derivative, so that the $f_{\bs{z}}$ is unimodal.\vs
\noindent
\tb{Transformation} ~
Assume the Cholesky decomposition of $\bs{\Sigma}$ yields a nonsingular lower triangular matrix $\bs{L}$ of the form
\[
\left(
\begin{array}{cc}
\bs{A} & \bs{0}\\
\bs{v}'& c
\end{array}
\right)
\]
with $\bs{\Sigma}=\bs{L}\bs{L}'$. Hence $\det(\bs{A})\neq 0\neq c$. Now transfer $(\bs{x}', y)$ to $(\bs{s}', t)$ with  $(\bs{s}', t)'=\bs{L}^{-1}((\bs{x}', y)'-\bs{\mu})$. It is readily seen that the distribution of $(\bs{s}', t)'$ follows
$E(g; \bs{0}, \bs{I_{p\times p}})$. \vs
Note that $(\bs{x}', y)'=\bs{L}(\bs{s}', t)'+(\bs{\mu}'_1, \mu_2)'$ with $\bs{\mu}=(\bs{\mu}'_1, \mu_2)'$. That is,
\begin{align}
\bs{x}&=\bs{A}\bs{s}+\bs{\mu}_1,\\[1ex]
y&=\bs{v}'\bs{s}+ct+\mu_2.
\end{align}
Equivalently, 
\begin{align}
(1,\bs{s}')'&=\bs{B}^{-1}(1,\bs{x}')',   \label{x-transformation.eqn}\\[1ex]
t&=\frac{y-(1,\bs{s}')(\mu_2, \bs{v}')'}{c}, \label{y-transformation.eqn}
\end{align}
where
\[
\bs{B}=\begin{pmatrix}
1 &\bs{0}'\\
\bs{\mu}_1  &\bs{A}
\end{pmatrix}
,~~~~
\bs{B}^{-1}=
\begin{pmatrix}
1 &\bs{0}'\\
-\bs{A}^{-1}\bs{\mu}_1& \bs{A}^{-1}
\end{pmatrix},
\]
\vs
It is readily seen that (\ref{x-transformation.eqn}) is an affine transformation on $\bs{w}$ and (\ref{y-transformation.eqn}) is first an affine transformation on $\bs{w}$ then a regression transformation on $y$ followed by a scale transformation on $y$. In light of Theorem 2.4, we can assume hereafter, w.l.o.g. that $(\bs{x}', y)$ follows an $E(g; \bs{0}, \bs{I}_{p\times p})$ (spherical) distribution and $\bs{I}_{p\times p}$ is the covariance matrix of $(\bs{x}', y)$.
\vs
\noindent
\tb{Theorem 4.1} Assume that
\bi
\item[(i)] assumptions of Theorem 2.3 hold;
\item[(ii)] $e\sim\mathcal{N}(0, \sigma^2)$ and $\bs{x}$ are independent.
\ei
Then
\bi
\item[(1)] $P\nabla=\bs{0}$ and $P(\nabla\nabla')=8\sigma^2 C \bs{I}_{p\times p}$,\\[1ex]
 with $C=-\alpha c \Phi'(\alpha c)+\Phi(\alpha c)-1/2$ where $\Phi$ is the cumultive distribution function of $\mathcal{N}(0,1)$ and $c=\Phi^{-1}(3/4)$.
\item[(2)] $\mb{V}= 2C_1\bs{I}_{p\times p} $ with $C_1=2*\Phi(\alpha c)-1$.
\item[(3)] $n^{1/2}(\widehat{\bs{\beta}}^n_{lst} - \bs{\beta}_{lst}) \stackrel{d} \longrightarrow  {\cal{N}}(O, \frac{2C\sigma^2}{C_1^2}\bs{I}_{p \times p})$
\ei
\vs
\noindent
\tb{Proof}: By Theorems 2.3 and 2.4, we can assume, w.l.o.g., that $\widehat{\bs{\beta}}_{lst}=\bs{\beta}_0=\bs{0}$. Utilizing the independence between $e$ and $\bs{x}$ and Theorem 3.2, a straightforward calculation leads to the results.
\hfill \pend
\vs
\noindent
\tb{Approximate $100(1-\gamma)\%$ confidence region} 
\vs
\noin
\tb{(i) Based on the asymptotic normality}
~Under the setting of Theorem 4.1, an approximate $100(1-\gamma)\%$ confidence region for the unknown regression parameter $\bs{\beta}_0$ is:
$$
\Big\{\bs{\beta} \in \R^p:~~ \|\bs{\beta}-\widehat{\bs{\beta}}^n_{lst}\|\leq \sqrt{\frac{2C\sigma^2}{C_1^2n}}\Phi^{-1}(\gamma)\Big \},
$$
where $\|\cdot\|$ stands for the Euclidean distance. Without the asymptotic normality, one can appeal to the next procedure.
\vs
\noin
\tb{(ii) Based on bootstrapping scheme and depth-median and depth-quantile}
This approximate procedure first re-samples $n$ points with replacement from the given original sample points and  calculates an $\widehat{\bs{\beta}}^n_{lst}$ (see Zuo (2022)) . Repreat this $m$ (a large number, say $10^4$) times and obtain $m$ such  $\widehat{\bs{\beta}}^n_{lst}$s. The next step is to calculate the depth, with respect to a location depth function ( e.g. halfspace depth (Zuo (2019)) or projection depth (Zuo (2003) and Shao and Zuo (2020)), of these $m$ points in the parameter space of $\bs{\beta}$.
Trimming $\lfloor\gamma m \rfloor$ of the least deepest points among the $m$ points, the left points form a convex hull, that is an approximate  $100(1-\gamma)\%$ confidence region for the unknown regression parameter $\bs{\beta}_0$ (see Zuo (2010, 2009) for the location case in low dimensions).

\section{Concluding remarks}

For  the establishment of the asymptotic normality (i.e. Theorem 3.2), the major contribution,
this article re-establishes some preliminary results in Section 2, some of those are established without any assumption (e.g. (i) of Theorems 2.1 and 2.2) and some with much less assumptions (e.g. Theorem 2.3 and (ii) of Theorem 2.2) and some established  with a different approach (e.g. Theorems 2.2 and 3.1). \vs
The asymptotic normality is applied in Theorem 4.1 for the practical inference procedure of confidence regions of the regression parameter $\bs{\beta}_0$. There are open problems left here, one is the estimation of the variance of $e$, which is now unrealistically assumed to be known, the other is the testing of hypothesis on $\bs{\beta}_0$.

\vs
\begin{center}
{\textbf{\large Acknowledgments}}
\end{center}

The author thanks  Prof. Wei Shao for
insightful comments and useful suggestions. 
\vs
\vs

\vs
\noin


\begin{thebibliography}{99}
{\small
\bibitem{B53} Box, G.E.P. (1953), ``Non-normality and tests on variances", {\it Biometrika}, 40, 318-335.


\bibitem{H64} Huber, P.\ J.\ (1964), ``Robust estimation of a location parameter'', \emph{Ann. Math. Statist.}, 35 73-101.
\bibitem{P84} Pollard, D. (1984), \emph{Convergence of Stochastic Processes}, Springer, Berlin.
\bibitem{R84} Rousseeuw,\ P.\ J. (1984), ``Least median of squares regression", {\it J. Amer. Statist. Assoc.} {\bf 79}, 871-880.
\bibitem{RH99} Rousseeuw, P.\ J., and Hubert, M. (1999),  ``Regression depth (with discussion)", \emph{J. Amer. Statist. Assoc.}, 94, 388--433.
\bibitem{RL87}Rousseeuw, P.J., and Leroy, A.  {\it Robust regression and outlier detection}. Wiley New York. (1987).
\bibitem{RY84} Rousseeuw, \ P.\ J.  and  Yohai, V. J. (1984), ``Robust regression by means of
S-estimators". In \emph{Robust and Nonlinear Time Series Analysis}. Lecture Notes
in Statist. Springer, New York. 26 256-272
\bibitem{SZ20} Shao, W. and Zuo, Y. (2020),
``Computing the halfspace depth with multiple try algorithm and simulated annealing algorithm", \emph{Comput Stat} 35, 203–226 (2020). https://doi.org/10.1007/s00180-019-00906-x
\bibitem{T60} Tukey, J. (1960), ``A survey on sampling from contaminated distributions", In \emph{Contributions to Probability and Statistics}. (I. Dlkin, ed.). 
Stanford University Press, Stanford, CA.
\bibitem{WZ09} Wu, M., and Zuo, Y. (2009), ``Trimmed and Winsorized means based on a scaled deviation",
\emph{J. Statist. Plann. Inference}, 139(2), 350-365.
\bibitem{Y87} Yohai, V.J. (1987), ``High breakdown-point and high efficiency estimates for regression", \emph{Ann. Statist.}, 15, 642–656.

\bibitem{YZ88} Yohai, V.J. and Zamar, R.H. (1988), ``High breakdown estimates of regression by means of
the minimization of an efficient scale", {\it J. Amer. Statist. Assoc.}, 83, 406–413.
\bibitem{Z03} Zuo, Y. (2003) ``Projection-based depth functions and associated medians'',
\emph{Ann. Statist.}, 31, 1460-1490.
\bibitem{Z06} Zuo, Y. (2006), ``Multi-dimensional trimming based on projection depth", \emph{Ann. Statist.}, 34(5), 2211-2251.
\bibitem{Z09} Zuo, Y. (2009), ``Data Depth Trimming Counterpart of the Classical $t$  (or $T^2$ ) Procedure", \emph{Journal of Probability and Statistics}, Volume 2009 |Article ID 373572 | https://doi.org/10.1155/2009/373572
\bibitem{Z10} Zuo, Y. (2010), ``Is the $t$ Confidence Interval $\overline{X}\pm t_{\alpha} (n-1)s/\sqrt{n}$ Optimal?'', \emph{The American Statistician}, 64:2, 170-173, DOI: 10.1198/tast.2010.09021
\bibitem{Z19}  Zuo, Y. (2019), ``A new approach for the computation of halfspace depth in high dimensions", \emph{Communications in Statistics - Simulation and Computation}, 48:3, 900-921, DOI: 10.1080/03610918.2017.1402040
\bibitem{Z20} Zuo, Y. (2020), ``Large sample properties of the regression depth induced median”,
\emph{Statistics and Probability Letters}, November 2020 166, arXiv1809.09896.
\bibitem{Z21a} Zuo, Y. (2021a), ``On general notions of depth for regression”
\emph{Statistical Science} 2021, Vol. 36, No. 1, 142–157, arXiv:1805.02046.

\bibitem{Z21b} Zuo, Y. (2021b), ``Robustness of the deepest projection regression depth functional",  {\it Statistical Papers}, vol. 62(3), pages 1167-1193.

\bibitem{Z22} Zuo, Y. (2022), ``Least sum of squares of trimmed residuals regression", arXiv:2202.10329
}
\end{thebibliography}

\begin{thebibliography}{99}
{\small

\bibitem{BC93} Bednarski, T. and Clarke, B.R. (1993), ``Trimmed likelihood estimation of location and scale of the normal
distribution". \emph{Austral. J. Statist.} 35, 141–153.



\bibitem{BV04} Boyd, S. and Vandenberghe, L. (2004),\emph{ Convex Optimization}. Cambridge University Press.

\bibitem{B82} Butler, R.W. (1982),  ``Nonparametric interval point prediction using data trimmed by a Grubbs type outlier
rule". \emph{Ann. Statist.} 10, 197–204.
\bibitem{DH83} Donoho,\ D.\ L., and Huber,\ P.\ J. (1983), ``The notion of
breakdown point", in: P.\ J.\ Bickel, K.\ A.\ Doksum and J.\ L.\ Hodges, Jr., eds. {\it A Festschrift
foe Erich L.\ Lehmann} (Wadsworth, Belmont, CA) pp.\ 157-184.

\bibitem{Eetal01} Edgar, T. F., Himmelblau, D. M., and Lasdon, L. S. (2001),  \emph{Optimization of Chemical Processes}, 2nd Edition,
McGraw-Hill Chemical Engineering Series.

\bibitem{H71} Hampel, F.R. (1971). ``A general qualitative definition of robustness", \emph{Ann. Math.
Statist.}, {{42}}, 1887-1896.

\bibitem{HRRS86} Hampel, F.\ R., Ronchetti, E.\ M., Rousseeuw, P.\ J., and Stahel, W.\ A. (1986),
\emph{Robust Statistics: The Approach Based on Influence Functions}, John Wiley \& Sons, New York.

\bibitem{H94} Hawkins, D. M.  (1994), ``The feasible solution algorithm
for least trimmed squares regression", \emph{Computational Statistics \& Data Analysis}, 17, 185-196.

\bibitem{HO99} Hawkins, D. M.  and Olive, D. J. (1999), "Improved feasible solution algorithms for high breakdown estimation "
\emph{Computational Statistics \& Data Analysis},   30(1), 1-11.

\bibitem{H67} Hodges, I.L. (1967). Efficiency in normal samples and tolerance of extreme values
for some estimates of location. Proceedings of the 5th Berkeley Symposium on
Mathematical Statistics and Probability, Vol. 1 pp. 163-168.

\bibitem{H94} H\"{o}ssjer, O. (1994), ``Rank-Based Estimates in the Linear Model with High Breakdown Point'',  {\it J. Amer. Statist. Assoc.}, {\bf 89}, 149-158.

\bibitem{H64} Huber, P.\ J.\ (1964), ``Robust estimation of a location parameter'', \emph{Ann. Math. Statist.}, 35 73-101.

\bibitem{L44} Levenberg, K. (1944), "A Method for the Solution of Certain Problems in Least Squares." \emph{Q. Appl. Math.}, 2, 164-168.
\bibitem{M63} Marquardt, D. (1963),  ``An Algorithm for Least-Squares Estimation of Nonlinear Parameters",  \emph{SIAM J. Appl. Math.}, 11, 431-441.

\bibitem{Netal89} Nierenberg DW, Stukel TA, Baron JA, Dain BJ, Greenberg ER(1989),  ``Determinants of plasma levels of beta-carotene and retinol",\emph{ American Journal of Epidemiology} 130:511-521.
\bibitem{Oetal15} \"{O}llerer, V., Croux, C., and Alfons, A. (2015) ``The
influence function of penalized regression estimators", \emph{Statistics}, 49:4, 741-765

\bibitem{P84} Pollard, D. (1984), \emph{Convergence of Stochastic Processes}, Springer, Berlin.
\bibitem{S87} Simpson, D.G. (1987),  
 ``Introduction to Rousseeuw (1984) Least Median of Squares Regression". In: Kotz S., Johnson N.L. (eds) \emph{Breakthroughs in Statistics}.  Springer Series in Statistics (Perspectives in Statistics). Springer, New York, NY. https://doi.org/10.1007/978-1-4612-0667-$5\_{18}$, 433-461
\bibitem{R83} Rousseeuw,\ P.\ J.  (1983), ``Multivariate estimation with high breakdown point". In \emph{Mathematical Statistics and
Applications}, eds W. Grossman, G. Pflug, I. Vincze, W. Wertz, pp. 283–297, Vol. B (1985). Dordrecht: Reidel Publishing Co.
\bibitem{R84} Rousseeuw,\ P.\ J. (1984), ``Least median of squares regression", {\it J. Amer. Statist. Assoc.} {\bf 79}, 871-880.
\bibitem{RH99} Rousseeuw, P.\ J., and Hubert, M. (1999),  ``Regression depth (with discussion)", \emph{J. Amer. Statist. Assoc.}, 94, 388--433.
\bibitem{RL87}Rousseeuw, P.J., and Leroy, A.  {\it Robust regression and outlier detection}. Wiley New York. (1987).
\bibitem{RD99}
Rousseeuw, P. J. and Van Driessen, K. (1999),
``A fast algorithm for the minimum covariance determinant estimator", \emph{Technometrics}, 41(3), 212-223.
\bibitem{RD06}
Rousseeuw, P. J. and Van Driessen, K. (2006), ``Computing LTS Regression for Large Data Sets", \emph{Data Mining and Knowledge Discovery} 12, 29-45.

\bibitem{RY84} Rousseeuw, \ P.\ J.  and  Yohai, V. J. (1984), ``Robust regression by means of
S-estimators". In \emph{Robust and Nonlinear Time Series Analysis}. Lecture Notes
in Statist. Springer, New York. 26 256-272

\bibitem{T94} Tableman, M. (1994),   ``The influence functions for the least
trimmed squares and the least trimmed absolute deviations estimators",  \emph{Statistics \& Probability Letters} 19 (1994) 329-337.
\bibitem{V06a}
V\'{ı}\v{s}ek, J.  \'{A}. (2006a), The least trimmed squares. Part I: Consistency.  \emph{Kybernetika},  42, 1-36.

\bibitem{V06b} V\'{ı}\v{s}ek, J.  \'{A}. (2006b) The least trimmed squares. Part II:
$\sqrt{n}$-consistency. \emph{Kybernetika},  42, 181-202.

\bibitem{V06c} V\'{ı}\v{s}ek, J.  \'{A}. (2006c), The least trimmed squares. Part III: Asymptotic normality.  \emph{Kybernetika}, 42,
203-224.

\bibitem{Y87} Yohai, V.J. (1987), ``High breakdown-point and high efficiency estimates for regression", \emph{Ann. Statist.}, 15, 642–656.

\bibitem{YZ88} Yohai, V.J. and Zamar, R.H. (1988), ``High breakdown estimates of regression by means of
the minimization of an efficient scale", {\it J. Amer. Statist. Assoc.}, 83, 406–413.

\bibitem{Z20} Zuo, Y. (2020), ``Large sample properties of the regression depth induced median”,
\emph{Statistics and Probability Letters}, November 2020 166, arXiv1809.09896.

\bibitem{Z21a} Zuo, Y. (2021a), ``On general notions of depth for regression”
\emph{Statistical Science} 2021, Vol. 36, No. 1, 142–157, arXiv:1805.02046.

\bibitem{Z21b} Zuo, Y. (2021b), ``Robustness of the deepest projection regression depth functional",  {\it Statistical Papers}, vol. 62(3), pages 1167-1193.

\bibitem{Z22} Zuo, Y. (2022), ``Least sum of squares of trimmed residuals regression", arXiv:2202.10329
}
\end{thebibliography}
\end{document}